\documentclass[a4,leqno,12pt]{article} 

\usepackage{amsmath,amsthm,amssymb,bm} 
\numberwithin{equation}{section}

\title{{\bf Regularized Product Expressions \\ 
of Higher Riemann Zeta Functions }} 
\author{Tetsuya MOMOTANI \\
\small
Graduate School of Mathematics, Kyushu University \\ 
\small 
6-10-1, Hakozaki Fukuoka 812-8581, Japan \\
\small 
E-mail: momo@math.kyushu-u.ac.jp
} 

\newtheorem{The}{Theorem}[subsection]
\newtheorem{Pro}[The]{Proposition}
\newtheorem{Lem}[The]{Lemma}
\newtheorem{Def}[The]{Definition}

\newtheorem{Cor}[The]{Corollary}

\newcommand{\minibullet}{%
\unitlength 0.1in 
\begin{picture}(0,0)(2.0,-2.4) 
\special{pn 12}%
\special{sh 1}%
\special{ar 200 200 12 12 0 6.2831853}%
\end{picture}}%

\DeclareMathOperator*{\rprod}    
{\displaystyle\coprod\kern-1.48em\prod} 
\DeclareMathOperator*{\grprod}   
{\displaystyle\coprod\kern-0.82em{\minibullet}\kern-0.82em\prod} 

\newcommand{\hprod}{\displaystyle\rprod} 
\newcommand{\mprod}{\displaystyle\grprod} 

\DeclareMathOperator*{\rrprod} 
{\coprod\kern-1.1em\prod} 
\DeclareMathOperator*{\ggprod} 
{\coprod\kern-0.66em{\minibullet}\kern-0.66em\prod}

\newcommand{\smprod}{\ggprod} 

\DeclareMathOperator*{\cterm}{CT} 

\begin{document} 

\maketitle 

\begin{abstract} 
As a generalization of \cite{KMW}, 
we introduce a higher Riemann zeta function 
for an abstract sequence. 
Then we explicitly determine 
its regularized product expression. 
\end{abstract} 

\section{Introduction} 

It is known (\cite{De}) that the Riemann zeta function 
$ \zeta (s) = \prod_{ p : \textit{prime} } (1 - p^{-s} )^{-1} $ 
has the regularized product expression: 
\begin{align*} 
\hprod_{ \rho \in R } 
\big( \frac{ s - \rho }{ 2 \pi } \big)
= 2^{-1/2} ( 2 \pi )^{-2} \pi ^{-s/2} 
s (s-1) \Gamma \big( \frac{s}{2} \big) \zeta (s) , 
\end{align*} 
where $ R $ is 
the set of the non-trivial zeros of $ \zeta (s) $, 
and $ \Gamma (s) $ is the classical gamma function. 

On the other hand,  
a higher Riemann zeta function 
$ \zeta_{ l \infty } (s) 
:= \prod_{ n=1 }^{ \infty } \zeta (s+ln) $ 
has been introduced and studied in the paper \cite{CL} for $l=1$, 
\cite{KMW} for $ l \in \mathbb{Z} _{\geq 1} $. 
We now consider a generalized higher Riemann zeta function. 
Let $ \Lambda = \{ \lambda_{k} \} _{ k \in I} $ 
be a sequence of complex numbers. 
Then we define 
a higher Riemann zeta function for the sequence $ \Lambda $ by 
\begin{align*} 
Z( s, \Lambda ) 
:= \prod_{ \lambda \in \Lambda } \prod_{ p : \textit{prime} }
( 1 - p^{ -s - \lambda } )^{-1} 
= \prod_{ \lambda \in \Lambda } \zeta ( s + \lambda ). 
\end{align*} 

In this paper, 
we study several properties of the higher Riemann zeta function. 
For a ``regularizable'' sequence $ \Lambda $, 
we see that $ Z(s, \Lambda) $ has the regularized product expression: 
\begin{align*} 
\mprod_{ \lambda \in \Lambda , \; \rho \in R } 
\big( \frac{ s + \lambda - \rho }{ 2 \pi } \big)
= \textit{`` gamma factor ''} \times Z( s, \Lambda ) . 
\end{align*}
Here $ \smprod $ denotes the ``dotted'' regularized product 
due to \cite{KW3}. 

Further, when $ \Lambda $ is given by a lattice 
$ \Omega 
:= \{ n_{1} \omega_{1} + \cdots + n_{r} \omega_{r} | 
n_{j} \in \mathbb{Z} _{ \geq 0} \} $, 
we show that the higher Riemann zeta function 
\begin{align*} 
Z ( s, \Omega ) 
= \prod _{ n_{1}, \cdots , n_{r} = 0 } ^{ \infty } 
\; \prod _{ p : \textit{prime} } 
( 1 - p^{ - ( s + n_{1} \omega_{1} + \cdots + n_{r} \omega_{r} ) } ) ^{-1} 
\end{align*} 
satisfies a certain functional equation of the type 
$ s \leftrightarrow 1 + \omega_{1} + \cdots + \omega_{r} -s $. 
We remark that 
the form of $ Z ( s, \Omega ) $ looks like the Selberg zeta function 
(when $r=1$)
defined for a discrete subgroup of 
a real rank $1$ semi-simple Lie group 
(Cf. \cite{Go} \cite{Ku2}). 

\section{Notation} 

\subsection{Regularized Product} 

Let 
$ \Lambda = \{ \lambda _{k} \} _{ k \in I } $ 
be a sequence of complex numbers. 
We call the sequence $ \Lambda $ is ``regularizable'' 
if $ \Lambda $ satisfies the following conditions: 
(i) 
$ \mathrm{Re} ( \lambda_{k} ) \geq 0 $ for all $ k \in I $, 
and $ \mathrm{Re} ( \lambda_{k} ) \to \infty $ as $ k \to \infty $. 
(ii) 
The series 
$ \sum_{ k \in I , \; \lambda_{k} \neq 0 } \lambda _{k} ^{-s} $ 
converges for sufficiently large $ \mathrm{Re} (s) $. 
(iii) 
The function 
$ \Theta ( x, \Lambda ) 
:= \sum _{ k \in I } e^{ - \lambda_{k} x }$ 
has an asymptotic expansion: 
\begin{align*} 
\Theta ( x, \Lambda ) 
\sim 
\sum_{ n=0 }^{ \infty } 
x^{ t_{n} } T_{n} ( \log x ) 
\quad \text{as } x \to + 0 , 
\end{align*} 
where $ t_{0} \leq 0 $, 
$ t_{0} < t_{1} < t_{2} < \cdots \to + \infty $, 
and $ T_{n} (z) \in \mathbb{C} [z] $ is a polynomial. 

For a regularizable sequence, 
we can define the regularized product as follows. 

\begin{Lem}[\cite{Ill}, \cite{KW3}] 
Assume that the sequence 
$ \Lambda = \{ \lambda _{k} \} _{ k \in I } $ is regularizable. 
Then, 
(1) 
for $ \mathrm{Re} (z) > 0 $, 
the function 
\begin{align*} 
\zeta ( s, z, \Lambda ) 
:= \sum_{ k \in I } ( z + \lambda _{k} ) ^{ -s } 
\end{align*} 
has an analytic continuation 
as a meromorphic function on $ s \in \mathbb{C} $. 
(2) The ``dotted'' regularized product 
\begin{align*} 
\mprod_{ k \in I } ( z + \lambda_{k} ) 
:= \exp \Big\{ 
- \cterm _{ s=0 } 
\frac{ \zeta ( s, z, \Lambda ) }{ s }
\Big\} 
\end{align*} 
exists, and this is an entire function on $ z \in \mathbb{C} $ 
with zeros at $ z = - \lambda_{k} $ $ ( k \in I ) $. 
Here $ \cterm _{s=0} f(s)$ denotes the constant term in 
the Laurent expansion of $ f(s) $ at $ s=0 $. 
\qed 
\end{Lem} 

We remark that the dotted regularized product 
$ \smprod _{ m \in I , \; n \in J } ( z + a_{m} + b_{n} ) $ 
exists 
if $ \{ a_{ m } \} _{ m \in I } $ and $ \{ b_{ n } \} _{ n \in J } $ 
are regularizable sequences. 
This follows at once from the relation 
$ \Theta ( x, \{ a_{m} + b_{n} \}_{ m \in I, \; n \in J } ) 
= \Theta ( x, \{ a_{ m } \} _{  m \in I } ) \times 
\Theta ( x, \{ b_{n} \} _{ n \in J } ) $. 

\subsection{Multiple Gamma Function and Multiple Sine Function} 

Let $ { \bm \omega } = ( \omega_{1}, \cdots , \omega_{r}) $ 
where $ \omega_{j} \in \mathbb{C} $, $ \mathrm{Re} ( \omega_{j} ) > 0 $, 
and 
$ \Omega := \{ n_{1} \omega_{1} + \cdots + n_{r} \omega_{r} | 
n_{j} \in \mathbb{Z} _{ \geq 0} \} $. 
We fix $ \mathrm{Re} (z) > 0 $. 
Then the multiple zeta function of Barnes (\cite{Ba}) is defined by 
\begin{align*} 
\zeta ( s, z, {\bm \omega} ) := \zeta ( s, z, \Omega )
= \sum_{ n_{1} , \cdots , n_{r} = 0 }^{ \infty } 
( z + n_{1} \omega_{1} + \cdots + n_{r} \omega_{r} )^{ -s } . 
\end{align*} 
This sum converges absolutely for $ \mathrm{Re} (s) > r $. 
It is seen that the sequence $ \Omega $ is regularizable, 
and $ \zeta ( s, z, { \bm \omega } ) $ 
has a meromorphic continuation to $ s \in \mathbb{C} $. 
Further, for $ m \in \mathbb{Z} _{ \geq 0} $, it is known that 
\begin{align*} 
\zeta ( -m , z, { \bm \omega } ) 
= \frac{ (-1)^{m} \cdot m! \cdot B_{m+r} ( z, { \bm \omega } ) } { (m+r)! }. 
\end{align*} 
Here $ B_{n} (z, { \bm \omega } ) $ 
is the multiple Bernoulli polynomials given by 
\begin{align*} 
\frac{ x^{r} e^{-zx} }
{ (1-e^{ - \omega_{1} x } ) \cdots ( 1-e^{ - \omega_{r} x } ) }
= \sum_{ n=0 }^{ \infty } 
\frac{ B_{n} (z, { \bm \omega } ) }{ n! } x^{n} , 
\end{align*} 
for $ |x| < \min_{1 \leq j \leq r} | 2\pi / \omega_{j} | $. 
We remark that 
$ \zeta (s) = \zeta ( s, 1, (1) ) $ is the Riemann zeta function and 
$ \zeta ( -m ) = (-1)^{m} B_{m+1} / ( m+1 ) $, 
where $ B_{m} = B_{m} (1,(1) ) $ is the usual Bernoulli number. 

Now we define 
the multiple gamma function $ \Gamma ( z, { \bm \omega } ) $ 
and the multiple sine function $ S ( z, { \bm \omega } ) $ 
by using the regularized product as follows (\cite{Ku1}). 
\begin{align*} 
\Gamma ( z, { \bm \omega } ) 
:=& \exp \Big\{ 
\frac{d}{ds} \zeta ( s, z, { \bm \omega } ) \Big| _{ s=0 } 
\Big\} 
= \hprod_{ n_{1} , \cdots , n_{r} =  0 }^{ \infty } 
( z + n_{1} \omega_{1} + \cdots + n_{r} \omega_{r} ) ^{ -1 } , \\
 S ( z, { \bm \omega } ) 
:=& \Gamma ( z , { \bm \omega } ) ^{-1} \cdot 
\Gamma ( \omega_{1} + \cdots + \omega _{r} - z, { \bm \omega} )^{(-1)^{r}} . 
\end{align*} 
It is seen that these functions satisfy the condition: 
\begin{align*} 
& \Gamma ( z, ( \omega_{1} , \cdots , \omega_{r-1} , \omega_{r} ) ) 
= \Gamma ( z, ( \omega_{1} , \cdots , \omega_{ r-1 } ) ) \cdot 
\Gamma ( z + \omega_{r} , 
( \omega_{1} , \cdots , \omega_{ r-1 } , \omega_{r} ) ) , \\ 
& S ( z, ( \omega_{1} , \cdots , \omega_{r-1} , \omega_{r} ) ) 
= S ( z, ( \omega_{1} , \cdots , \omega_{ r-1 } ) ) \cdot 
S ( z + \omega_{r} , ( \omega_{1} , \cdots , \omega_{ r-1 } , \omega_{r} ) ) . 
\end{align*} 
And we note that $ \Gamma (z, (1) ) = ( 2 \pi ) ^{-1/2} \; \Gamma (z) $ and 
$ S ( z, (1) ) = 2 \sin ( \pi z ) $. 

\section{Higher Riemann Zeta Function} 

\subsection{Definition and Analytic Continuation} 

First of all, we introduce a higher Riemann zeta function as follows. 

\begin{Def} 
Define a higher Riemann zeta function for 
$ \Lambda = \{ \lambda_{k} \} _{ k=0 }^{ \infty } $ by 
\begin{align*} 
Z ( s , \Lambda ) 
= Z ( s , \{ \lambda_{k} \} _{ k=0 }^{ \infty } )
:= \prod _{ k=0 }^{ \infty } \zeta ( s + \lambda _{k} ) 
= \prod _{ k=0 }^{ \infty } \prod_{ p : \text{prime} } 
( 1 - p^{ - ( s + \lambda_{k} ) } ) ^{-1} .
\end{align*} 
\end{Def} 

To show the absolute convergence of $ Z ( s , \Lambda ) $, 
we prepare the following lemma. 

\begin{Lem} 
Fix a positive number $ x \geq 2 $. 
If the sequence $ \Lambda = \{ \lambda_{k} \} _{ k=0 }^{ \infty } $ satisfies 
(i) $ 0 \leq \mathrm{Re} ( \lambda _{0} ) 
\leq \mathrm{Re} ( \lambda_{1} ) \leq \cdots \to \infty $, 
and (ii) the sum 
$ \sum_{ k \in I , \; \lambda_{k} \neq 0 } \lambda_{k} ^{-s} $ 
converges absolutely for large $ \mathrm{Re}(s) $, 
then we have 
\begin{align*}
\big| \sum_{ k=j }^{ \infty } x ^{ - \lambda_{k} } \big| 
< \hspace{-1mm} < _{j} x ^{ - \mathrm{Re} ( \lambda_{j} ) }, 
\qquad 
\text{ for all $j \geq 0$ }. 
\end{align*} 
\end{Lem} 

\begin{proof} 
We observe that 
\begin{align*} 
\big| \sum_{ k=j }^{ \infty } x ^{ - \lambda_{k} } \big|
=& \big| x ^{ - \lambda_{j} } 
\cdot \sum_{ k=j }^{ \infty } x^{ - ( \lambda_{k} - \lambda_{j} ) } \big| \\ 
\leq & x ^{ - \mathrm{Re} ( \lambda_{j} ) } 
\cdot \sum_{ k=j }^{ \infty } 
2 ^{ - ( \mathrm{Re} ( \lambda_{k} )- \mathrm{Re} ( \lambda_{j} ) ) } .
\end{align*} 
Under the condition (ii), the sum of the second factor is bounded. 
\end{proof} 

\begin{The} 
Assume that the sequence 
$ \Lambda = \{ \lambda_{k} \} _{ k=0 }^{ \infty } $ satisfies 
the conditions (i) and (ii) in previous lemma. 
Then $ Z ( s , \Lambda ) $ is defined for 
$ \mathrm{Re} (s) > 1 - \mathrm{Re} ( \lambda _{0} ) $. 
Moreover $ Z ( s, \Lambda ) $ has an analytic continuation 
as a meromorphic function on $ s \in \mathbb{C} $. 
\end{The} 

\begin{proof}
We see from the lemma that 
\begin{align*} 
\big| \sum_{ k=0 }^{ \infty } \sum_{ p }  
p ^{ - ( s + \lambda _{k}) } \big| 
\ll \sum_{ p } p ^{ - \mathrm{Re} ( s ) -\mathrm{Re} ( \lambda _{0} ) } .
\end{align*} 
The right side converges in 
$ \mathrm{Re} (s) > 1 - \mathrm{Re} ( \lambda _{0} ) $. 
Hence $ Z ( s , \Lambda ) $ is holomorphic on this domain. 
Further, we have 
\begin{align*} 
Z ( s , \{ \lambda _{k} \} _{ k=0 }^{ \infty } ) 
= \prod_{ k = 0 }^{ j-1 } \zeta ( s + \lambda _{k} ) 
\cdot \prod_{ k = j }^{ \infty } \zeta ( s + \lambda _{k} ) 
= \prod_{ k = 0 }^{ j-1 } \zeta ( s + \lambda _{k} ) 
\times Z ( s, \{ \lambda _{k} \} _{ k=j }^{ \infty } )
\end{align*} 
for any $j$. 
From the previous lemma, 
we see that the second factor is defined now for 
$ \mathrm{Re} (s) > 1 - \mathrm{Re} ( \lambda _{j} ) $. 
This provides a meromorphic continuation to $ \mathbb{C} $. 
\end{proof} 

We remark that 
the higher Riemann zeta function has the Dirichlet series expression 
$ Z ( s , \Lambda ) 
= \sum_{ n=1 } ^{ \infty } g_{ \Lambda } (n) n^{ -s } $ where 
\begin{align*} 
g _{ \Lambda } ( n ) 
= \sum_{ 
\substack { n_{0} \cdot n_{1} \cdot n_{2} \cdot \cdots \; = \, n \\ 
n_{0} \geq 1, \; n_{1} \geq 1, \cdots 
} } n_{0}^{ - \lambda _{0} } 
\cdot n_{1} ^{ - \lambda _{1} } 
\cdot n_{2} ^{ - \lambda _{2} } \cdots . 
\end{align*} 
From this expression, 
we see that 
$ g_{ \Lambda }(n) $ is multiplicative and 
$ Z (s, \Lambda ) $ has the Euler product: 
\begin{align*} 
Z ( s, \Lambda )  = \prod_{ p : \textit{prime} }
\sum_{ m=0 }^{ \infty } \frac{ g_{ \Lambda } (p^{m}) }{ p^{ms} } , \qquad  
\textit{ for } \mathrm{Re} (s) > 1 - \mathrm{Re} ( \lambda_{0} ) . 
\end{align*} 
Moreover, using the Tauberian theorem (Cf.\cite{Mu}), 
we obtain the behavior of $ g_{ \Lambda } (n) $ as follows. 

\begin{Cor} 
If $ \Lambda = \{ \lambda _{k} \} _{k=0} ^{\infty} $ satisfies 
(i) (ii) and 
$ 0 \leq \lambda_{0} = \cdots = \lambda_{K-1} < \lambda_{K} \leq \cdots $. 
Then we have 
\begin{align*} 
\sum_{ n \leq x } g_{ \Lambda } (n) 
= \big{(} c_{ \Lambda } + o (1) \big{)} 
x^{ 1 - \lambda_{0} } \cdot \log ^{ K -1 } x  \qquad \text{as } x \to \infty ,
\end{align*} 
where $ c_{ \Lambda } $ is a constant which is given by 
\begin{align*} 
c_{ \Lambda } := &  \frac{1}{ ( K -1 ) ! } 
\lim_{ s \to 1 - \lambda _{0} } ( s-1+ \lambda _{0} )^{ K } Z(s, \Lambda) \\
= & \frac{ 1 } { ( K -1 )! } 
Z( 1 - \lambda _{0 } , \{ \lambda _{k} \} _{ k=K }^{ \infty } ) .
\end{align*} 
\qed 
\end{Cor} 

\subsection{Regularized Product Expression} 

We first recall 
the regularized product expression of the Riemann zeta function. 

\begin{Lem}[\cite{De}] 
For $ \mathrm{Re} (z) > 1 $ and $ \mathrm{Re} (s) > 1 $, we have 
\begin{align} 
\sum_{ \rho \in R } ( z - \rho )^{ -s } 
= z^{-s} + (z-1)^{-s} - \sum_{ n=0 }^{ \infty } ( z+2n ) ^{-s} 
- \frac{ 1 }{ \Gamma (s) } 
\sum_{ k=0 }^{ \infty } \sum_{ p : \, \textit{prime} } 
\frac{ \log p }{ p^{ nz } } ( \log p^{n} )^{s-1} , 
\label{AAA} 
\end{align} 
where $ R $ denotes 
the set of the non-trivial zeros of the Riemann zeta function. 
\qed 
\end{Lem} 
This relation follows by applying Weil's explicit formula. 
For a proof, see \cite{De}. 
From the lemma we see that 
the function $ \sum_{ \rho \in R } ( z - \rho ) ^{-s} $ 
has a meromorphic continuation to $ s \in \mathbb{C} $ and 
the Riemann zeta function has the regularized product expression: 
\begin{align*} 
\hprod_{ \rho \in R } 
\big( \frac{ s - \rho }{ 2\pi } \big) 
= \frac{ s }{ 2 \pi } \cdot \frac{ s -1}{ 2 \pi } 
\cdot \hprod_{ n=0 }^{ \infty } \big( \frac{ s +2n }{ 2 \pi } \big) ^{-1} 
\cdot \zeta (s).
\end{align*} 
Further using properties of the Hurwitz zeta function 
and the gamma function, 
we have 
\begin{align*} 
\hprod_{ \rho \in R } \big( \frac{ s - \rho }{ 2\pi } \big) 
= 2^{-1/2} \; ( 2 \pi ) ^{ -2 } \cdot
s (s-1) \cdot \pi ^{ -s/2 } \cdot \Gamma ( \frac{s}{2} ) \cdot \zeta (s) . 
\end{align*} 
We remark that this function is invariant under 
$ s \leftrightarrow 1 - s $. 

Next we describe the regularized product expression 
of the higher Riemann zeta function as follows. 

\begin{The} 
Let 
$ \Lambda $ be a regularizable sequence. 
Then the higher Riemann zeta function $ Z (s , \Lambda ) $ 
has the following 
regularized product expression. 
\begin{align} 
\mprod_{ \lambda \in \Lambda , \; \rho \in R } 
\big( \frac{ s + \lambda - \rho }{ 2\pi } \big) 
= \mprod_{ \lambda \in \Lambda } 
\big( \frac{ s + \lambda }{ 2 \pi } \big) 
\cdot \mprod_{ \lambda \in \Lambda } 
\big( \frac{ s -1 + \lambda }{ 2 \pi } \big) 
\cdot \mprod_{ \lambda \in \Lambda , \; n \geq 0 } 
\big( \frac{ s + \lambda +2n }{ 2 \pi } \big) ^{-1} 
\cdot Z (s, \Lambda ). 
\label{BBB} 
\end{align} 
\end{The} 

\begin{proof} 
From (\ref{AAA}), we see that 
\begin{align*} 
\sum _{ \lambda \in \Lambda } \sum_{ \rho \in R } 
( z + \lambda - \rho )^{ -s } 
=& \sum _{ \lambda \in \Lambda } ( z + \lambda) ^{ -s } 
+ \sum _{ \lambda \in \Lambda } ( z + \lambda -1) ^{ -s } 
- \sum _{ \lambda \in \Lambda } \sum_{ n=0 }^{ \infty } 
( z + \lambda+ 2n ) ^{-s} \\ 
-& \frac{ 1 }{ \Gamma (s) } 
\sum _{ \lambda \in \Lambda } \sum_{ n=1 }^{ \infty } \sum_{ p } 
\frac{ \log p }{ p^{ n (z+\lambda) } } ( \log p^{n} )^{s-1} , 
\end{align*} 
for $ \mathrm{Re} (z) > 1 $ and large $ \mathrm{Re} (s) $. 
We observe that 
the right hand side is now meromorphic for any $ s \in \mathbb{C} $. 
Therefore the dotted product 
$ \smprod _{ \lambda , \; \rho } \{ ( z + \lambda - \rho ) / 2 \pi \} $ 
exists. 
Further we obtain the equation (\ref{BBB}) by using the relation 
\begin{align*} 
\log Z ( z , \Lambda ) 
= \sum_{ \lambda \in \Lambda } \sum_{ p } \sum_{ n=1 }^{ \infty } 
\frac{1}{ n p^{ n ( z + \lambda ) } } 
= \cterm _{ s= 0 } \frac{1}{s} \Big\{ 
\frac{ ( 2 \pi )^{s} }{ \Gamma (s) } 
\sum _{ \lambda \in \Lambda } \sum_{ n=1 }^{ \infty } \sum_{ p } 
\frac{ \log p }{ p^{ n (z+\lambda) } } ( \log p^{n} )^{s-1} \Big\} . 
\end{align*} 
This completes the proof. 
\end{proof} 

\subsection{Semi-Lattice and Functional Equation} 

Let 
$ \Omega := \{ n_{1} \omega_{1} + \cdots + n_{r} \omega_{r} | 
n_{j} \in \mathbb{Z} _{ \geq 0} \} $ 
with $ \mathrm{Re} ( \omega_{j} ) > 0 $. 
Then it is seen that $ \Omega $ is a regularizable sequence. 
Now, we consider the higher Riemann zeta function 
for the semi-lattice $ \Omega $. 

\begin{Def} 
Define the higher Riemann zeta function of the weight 
$ { \bm \omega } = ( \omega_{1} \cdots \omega_{r} ) $ 
by 
\begin{align*} 
Z ( s, { \bm \omega } ) 
:= Z ( s , \Omega ) 
= \prod _{ n_{1}, \cdots , n_{r} = 0 } ^{ \infty } 
\prod_{ p : \text{prime} } 
( 1 - p^{ - ( s + n_{1} \omega_{1} + \cdots + n_{r} \omega_{r} ) } ) ^{-1} . 
\end{align*} 
\end{Def} 
We see that this product absolutely for $ \mathrm{Re} (s) > 1 $ 
and $ Z ( s, { \bm \omega } ) $ has an analytic continuation 
as a meromorphic function on $ s \in \mathbb{C} $. 
Obviously, this is invariant under the arrangement of $ \omega_{j} $'s. 

\begin{Pro} 
The higher Riemann zeta function $ Z( s, { \bm \omega } ) $ 
defined for $ Re (s) > 1 $ 
has a meromorphic continuation to the whole complex plane 
and satisfies 
\begin{align} 
& Z ( s, ( \omega_{1} , \cdots , \omega_{ r-1 } , \omega_{r} ) ) 
= Z ( s, ( \omega_{1} , \cdots , \omega_{ r-1 } ) ) \cdot 
Z ( s + \omega_{r} , ( \omega_{1}, \cdots , \omega_{r-1}, \omega_{r} ) ) , 
\label{CCC} 
\\ 
& Z (s , { \bm \omega } ) \times 
\prod_{ 
\substack{ 1 \leq k_{1} < \cdots < k_{j} \leq r \\ 
1 \leq j \leq r } }  
Z (s + \omega_{ k_{1} } + \cdots + \omega_{ k_{j} } , 
{ \bm \omega } ) ^{ (-1)^{j} } = \zeta (s) . 
\notag 
\end{align} 
\end{Pro} 

\begin{proof} 
We observe that 
\begin{align*} 
Z( s, ( \omega_{1} , \cdots , \omega_{r-1} , \omega_{r} ) ) 
& = \prod_{ n=0 }^{ \infty } 
Z(s+ n \omega_{r} , ( \omega_{1} , \cdots , \omega_{r-1} )) \\ 
& = Z( s, ( \omega_{1} , \cdots , \omega_{ r-1 } ) ) 
\cdot \sum _{ n = 0 }^{ \infty } 
Z(s+ n \omega_{r} + \omega_{r} , ( \omega_{1}, \cdots , \omega_{r-1} )) \\ 
& = Z( s, ( \omega_{1} , \cdots , \omega_{ r-1 } ) ) \cdot 
Z( s + \omega_{r} , ( \omega_{1}, \cdots , \omega_{r-1}, \omega_{r} ) ) . 
\end{align*} 
Thus the relation (\ref{CCC}) holds. 
By using (\ref{CCC}) repeatedly, we have 
\begin{align*} 
\zeta (s) 
=& Z ( s , ( \omega_{1} ) ) \cdot Z ( s + \omega_{1} , ( \omega_{1} ) ) ^{-1} \\
=& Z ( s , ( \omega_{1} ,\omega_{2}) ) 
Z ( s + \omega_{2} , ( \omega_{1}, \omega_{2} ) ) ^{-1}
\cdot Z ( s + \omega_{1} , ( \omega_{1}, \omega_{2} ) ) ^{-1} 
Z ( s + \omega_{1}+ \omega_{2} , ( \omega_{1}, \omega_{2} ) ) \\
= & \cdots \\ 
= & Z (s , { \bm \omega } ) \times 
\prod_{ 
\substack{ 1 \leq k_{1} < \cdots < k_{j} \leq r \\ 
1 \leq j \leq r } }  
Z (s + \omega_{ k_{1} } + \cdots + \omega_{ k_{j} }, 
{ \bm \omega } ) ^{ (-1)^{j} } . 
\end{align*} 
This completes the proof. 
\end{proof} 

Next we show the regularized product expression of 
$ Z( s, { \bm \omega } ) $
as follows. 

\begin{The} 
We have 
\begin{align*} 
\hprod_{ n_{1} , \cdots , n_{r} \geq 0 , \; \rho \in R } & 
\Big( 
\frac{ s + n_{1} \omega_{1} + \cdots + n_{r} \omega_{r} - \rho } 
{ 2 \pi } \Big) \\ 
=& \exp \Big\{ \Big( 
- \frac{ B_{r} ( s , { \bm \omega } )}{ r! } 
- \frac{ B_{r} ( s-1 , { \bm \omega } )}{ r! } 
+ \frac{ B_{r+1} ( s , ( 2, \omega_{1} , \cdots , \omega_{r} ) ) }{ (r+1)! } 
\Big) \log ( 2 \pi ) \Big\} \\ 
\times & 
\Gamma ( s , { \bm \omega } ) ^{-1} 
\cdot \Gamma ( s-1 , { \bm \omega } ) ^{-1} 
\cdot \Gamma ( s , ( 2, \omega_{1} , \cdots , \omega_{r} ) ) 
\cdot Z ( s, { \bm \omega } ) . 
\end{align*} 
Here $ B_{n} (s, { \bm \omega } ) $ is the multiple Bernoulli polynomial, 
$ \Gamma ( s ,{ \bm \omega } ) $ is the multiple gamma function. 
\end{The} 

\begin{proof} 
We note that 
\begin{align*} 
\hprod_{ n_{1} , \cdots , n_{r} = 0 }^{ \infty } 
\Big( 
\frac{ z + n_{1} \omega_{1} + \cdots + n_{r} \omega_{r} }{ 2 \pi } 
\Big) 
= ( 2 \pi )^{ - \zeta ( 0, z, { \bm \omega } ) }
\hprod_{ n_{1} , \cdots , n_{r} = 0 }^{ \infty } 
( z + n_{1} \omega_{1} + \cdots + n_{r} \omega_{r} ) . 
\end{align*} 
Then, the result follows from Theorem 3.2.2. 
\end{proof} 

From this theorem, we see that the function 
\begin{align*} 
\hat{Z} (s , { \bm \omega } ) 
:= & \exp \Big\{ \Big( 
- \frac{ B_{r} ( s , { \bm \omega } )}{ r! } 
- \frac{ B_{r} ( s-1 , { \bm \omega } )}{ r! } 
+ \frac{ B_{r+1} ( s , ( 2, \omega_{1} , \cdots , \omega_{r} ) ) }{ (r+1)! } 
\Big) \log ( 2 \pi ) \Big\} \\ 
\times & 
\Gamma ( s , { \bm \omega } ) ^{-1} 
\cdot \Gamma ( s-1 , { \bm \omega } ) ^{-1} 
\cdot \Gamma ( s , ( 2, \omega_{1} , \cdots , \omega_{r} ) ) 
\cdot Z ( s, { \bm \omega } ) 
\end{align*} 
is an entire function of order $r+1$ with zeros at 
$ s = \rho - n_{1} \omega_{1} - \cdots - n_{r} \omega_{r} $. 
Further $ \hat{Z} (s , { \bm \omega } ) $ satisfies the relation 
\begin{align*} 
\hat{Z} ( s, ( \omega_{1} , \cdots , \omega_{ r-1 } , \omega_{r} ) ) 
= \hat{Z} ( s, ( \omega_{1} , \cdots , \omega_{ r-1 } ) ) \cdot 
\hat{Z} ( s + \omega_{r} , ( \omega_{1}, \cdots , \omega_{r-1}, \omega_{r} ) ) .
\end{align*}

Finally, 
we give the functional equation of the higher Riemann zeta function 
$ Z (s , {\bm \omega} ) $ as follows. 

\begin{The} 
Define the function $ \Lambda ( s, { \bm \omega } ) $ by 
\begin{align*}  
\Lambda ( s, { \bm \omega } ) 
:= \hat{Z} ( s , { \bm \omega } ) 
\cdot \hat{Z} 
( 1 + \omega_{1} + \cdots + \omega_{r} - s , { \bm \omega } )
^{ (-1)^{r+1} } . 
\end{align*}
Then we have 
\begin{align} 
\Lambda ( s, ( \omega_{1} , \cdots , \omega_{r-1} , \omega_{r} ) ) 
= \Lambda ( s, ( \omega_{1} , \cdots , \omega_{ r-1 } ) ) \cdot 
\Lambda 
( s + \omega_{r} , ( \omega_{1}, \cdots , \omega_{r-1}, \omega_{r} ) ) . 
\label{DDD} 
\end{align} 
Moreover, $ \Lambda (s , { \bm \omega } ) $ satisfies the functional equation: 
\begin{align} 
\Lambda (s , { \bm \omega } ) \times 
\prod_{ 
\substack{ 1 \leq k_{1} < \cdots < k_{j} \leq r \\ 
1 \leq j \leq r } }  
\Lambda (s + \omega_{ k_{1} } + \cdots + \omega_{ k_{j} }, 
{ \bm \omega } ) ^{ (-1)^{j} } 
= 1 . 
\end{align} 
\end{The} 

\begin{proof} 
we observe (\ref{CCC}) and 
\begin{align*}
& \hat{Z} ( 1 + \omega_{1} + \cdots + \omega_{r-1} + \omega_{r} - s , 
( \omega_{1} , \cdots , + \omega_{r-1} , \omega_{r} ) ) \\
=& \frac{ \hat{Z} 
( 1 + \omega_{1} + \cdots + \omega_{r-1} + \omega_{r} - ( s + \omega_{r} ) ,
( \omega_{1} , \cdots , + \omega_{r-1} , \omega_{r} ) ) }
{ \hat{Z} ( 1 + \omega_{1} + \cdots + \omega_{r-1} - s , 
( \omega_{1} , \cdots , \omega_{r-1} ) ) } .
\end{align*} 
Then the equation (\ref{DDD}) follows immediately. 

Next, using (\ref{DDD}) repeatedly, we have 
\begin{align*} 
\Lambda (s , { \bm \omega } ) \times 
\prod_{ 
\substack{ 1 \leq k_{1} < \cdots < k_{j} \leq r \\ 1 \leq j \leq r } }  
\Lambda (s + \omega_{ k_{1} } + \cdots + \omega_{ k_{j} }, 
{ \bm \omega } ) ^{ (-1)^{j} } 
= \hat{ \zeta } (s) \cdot \hat{ \zeta } (1-s) ^{-1} = 1 ,  
\end{align*} 
where 
\begin{align*}
\hat{ \zeta } (s)  
:= \hprod_{ \rho \in R } \big( \frac{ s - \rho }{ 2\pi } \big) 
= 2^{-1/2} \; ( 2 \pi ) ^{ -2 } \cdot
s (s-1) \cdot \pi ^{ -s/2 } \cdot \Gamma ( \frac{s}{2} ) \cdot \zeta (s) . 
\end{align*} 
Hence the theorem follows. 
\end{proof} 


\end{document}